\documentclass{amsart}   	
\usepackage{geometry}                		
\usepackage{graphicx,graphics, amsthm, amsfonts, amsbsy,amssymb, amsmath, cite, array, hyperref, enumerate}			
\renewenvironment{proof}{{\noindent\bfseries Proof.}}{\qed}								
\usepackage{amssymb}
\newtheorem{theorem}{Theorem}[section]

\newtheorem{example}[theorem]{Example}

\newtheorem{lemma}[theorem]{Lemma}
\newtheorem{corollary}[theorem]{Corollary}

\newtheorem{proposition}[theorem]{Proposition}

\author{Ralf Fr\"oberg}
\email{frobergralf@gmail.com}
\address{Department of Mathematics, Stockholm University, S-10691 Stockholm, Sweden}
\usepackage{xcolor}

\title{Some new Betti numbers of ideals generated by $n+1$ generic forms in $n$ variables}

\begin{document}

\begin{abstract}
Very little is known on the Hilbert series of graded algebras $\mathbb C[x_1,\ldots,x_n]/(g_1,\ldots,g_r)$, $r>n$,  $g_i$ generic form of degree $e_i$, 
in general. One instance when the series is known, is for $n+1$ forms in $n$ variables, \cite{St}. Of course even less is known about Betti numbers. There
are some general results on the Betti table by Pardue and Richert in \cite{Pa-Ri,Pa-Ri1}, and by Diem in \cite{Di}.
Then there are results on Betti numbers in the case $n+1$ relations in $n$ variables, described below, by Migliore and Mir\'{o}-Roig in \cite{Mi-Mi}.
In this paper we consider the same case as in \cite{Mi-Mi}, $n+1$ 
forms in $n$ variables. Our results can be described as follows. We can determine all graded Betti numbers of 
$\mathbb C[x_1,\ldots,x_n]/(g_1,\ldots,g_{n+1})$, $g_i$ generic, at least if $\sum_{i=1}^{n+1}\deg(g_i)-n$ is even, often in more cases, see Theorem \ref{par}. Thus,
given {\em any} set $\{ e_1,\ldots,e_n\}$, $e_i\ge2$ for all $i$, such that $\deg(g_i)=e_i$, $i=1,\ldots,n$, 
we get many numbers $e_{n+1}$, so that we can determine all graded Betti numbers of $\mathbb C[x_1,\ldots,x_n]/(g_1,\ldots,g_{n+1})$, 
$\deg(g_i)=e_i$, $1\le i\le n+1$. When we say that the Betti numbers could be determined, we don't claim that we can give formulas for them. 
We mean that the Betti numbers can be calculated from
the Hilbert series, which is known, and the Koszul complex on $g_1,\ldots,g_{n+1}$, which is well known and easy to get hold of. 
The main ingredients of the proof is a theorem by Pardue and Richert, \cite{Pa-Ri,Pa-Ri1},
and later by Diem,\cite{Di} together with a theorem on Hilbert function of artinian complete intersections by Reid, Roberts, and Roitman, \cite{R-R-R}.
A second goal of this article is to give a short elementary proof of the theorem by Reid, Roberts, and Roitman, see Theorem \ref{RRR}.
\end{abstract}

\maketitle

\section{Introduction}
This article is about Betti numbers of ideals generated by generic forms in polynomial rings. To avoid any difficulties we consider polynomial
rings over the complex numbers. In general not even the Hilbert series of ideals generated by generic forms is known. There is a conjecture:
If the ideal is $I=(g_1,\ldots,g_k)$ in $S=\mathbb C[x_1,\ldots,x_n]$, then the Hilbert series of $S/I$ is $[\prod_{i=1}^k(1-t^{e_i})/(1-t)^n]$, where
$e_i=\deg g_i$, and the brackets mean that we should only keep the terms as long as they have positive coefficients, \cite{Fr}.
The formula is trivially true if $k\le n$, since in that case we have a complete intersection, and it is proved to be true e.g.
for $n=3$ by Anick, \cite{An}, and for $k=n+1$ by Stanley, \cite{St}. There are general results on the minimal resolution of $S/I$ in \cite{Pa-Ri,Pa-Ri1}, \cite{Di},
and in \cite{Mi-Mi1}, and also determination of all graded Betti numbers in some cases when $r=n+1$ in \cite{Mi-Mi}. We shortly describe 
the results on Betti numbers below.

\section{What is known?}
In the sequel we assume $R$ to be artinian. For a graded module $M$, let $M_{(d)}$ be the submodule of $M$ generated by 
the elements of $M$ of degree at  most $d$.
Suppose that $C_*:\cdots\rightarrow C_1\rightarrow C_0$ is a complex and $e$ an integer such that for all $i$, 
$\delta_i^C((C_i)_{(e+i)})\subseteq(C_{i-1})_{(e+i-1)}$. Here $C_*$ could be the minimal resolution $F_*$ of $S/I$ or the Koszul complex $K_*$ on
$f_1,\ldots,f_k$. Then we have the restricted complex $C_*^{(e)}:\cdots\rightarrow(C_i)_{(e+i)}\rightarrow(C_{i-1})_{(e+i-1)}\rightarrow\cdots$.
Suppose $R=S/(f_1,\ldots,f_k)$, $f_i$ generic of degree $d_i$, with $R_d\ne0$ and $R_i=0$ if $i>d$, 
so $d$ is the highest socle degree of $R$, and also the regularity of $R$.
Then it is proved in \cite{Pa-Ri,Pa-Ri1} and \cite{Di} that $F_*^{(e)}$ is isomorphic to $K_*^{(e)}$ for $e\le d-2$. Furthermore, if the first nonpositive
coefficient in $\prod_{i=1}^k(1-t^{d_i})/(1-t)^n$ is 0, then $F_*^{(e)}$ is isomorphic to $K_*^{(e)}$ for $e\le d-1$.

\bigskip
Migliore and Mir\'{o}-Roig \cite{Mi-Mi} discuss the resolution of $\mathbb C[x_1,\ldots,x_n]/(g_1,\ldots,g_{n+1})$, where $g_i$ is generic of degree $e_i$. Among other things they 
determine all Betti numbers in the following six cases with $n+1$ generic relations in $n$ variables:
(1) $n=3$ (so 4 relations), (2) $n=4$ and $\sum_{i=1}^5e_i$ even, (3) $n=4$ and $\sum_{i=1}^5e_i$ odd and $e_2+e_3+e_4<e_1+e_5+4$,
(4) $n$ even and all generators of the same degree $a$ which is even, (5) $\sum_{i=1}^{n+1}e_i-n$ even and $e_2+\cdots+e_n<e_1+e_{n+1}+n$,
(6) $\sum_{i=1}^{n+1}e_i-n$ odd, $n\ge6$ even, $e_2+\cdots+e_n<e_1+e_{n+1}+n$ and $e_1+\cdots+e_n-e_{n+1}-n>>0$.

\medskip
Here is a (very coarse) description of their method. Let $S=k[x_1,\ldots,x_n]$, $I=(g_1,\ldots,g_{n+1})$ and $J=(g_1,\ldots,g_n)$ (a complete intersection).
The $I$ is linked to $G$ via $J$, i.e. $J:I=G$ and $J:G=I$, where $S/G$ is Gorenstein. Knowing the Hilbert functions of $S/I$ and $S/J$
one can determine the Hilbert function of $S/G$. The resolution of $S/G$ determines the
resolution of $S/I$. Suppose $S/G$ is compressed, i.e., that the Hilbert function is maximal given the socle degree and the
embedding dimension. Then the Hilbert function of $S/G$ is increasing to the middle and has one or two peaks. It has one peak if and only if
$\sum_{i=1}^{n+1}e_i-n$ (the socle degree) is even. $S/G$ is compressed if and only if $e_2+\cdots+e_n<e_1+e_{n+1}+n$. If it has one peak,
the resolution is pure, so the Betti numbers follow from the Hilbert series.

We repeat that when we say below that the Betti numbers could be determined, we don't claim that we can give formulas for them. We mean that the Betti numbers can be calculated from
the Hilbert series, which is known, and the Koszul complex on $g_1,\ldots,g_{n+1}$, which is well known and easy to get hold of. The reason for us not being able
to get formulas for the Betti numbers, is that it seems very hard to get formulas for the coefficients of $[\prod_{i=1}^{n+1}(1-t^{e_i})/(1-t)^n]$ in general.

\section{A new proof of the theorem by  Reid, Roberts, and Roitman.}
A sequence of numbers $a_0,a_1,\ldots,a_d$ is called unimodal if $a_0\le a_1\le\cdots\le a_k\ge a_{k+1}\ge\cdots\ge a_d$ for some $k$. It is called
symmetric if $a_i=a_{d-i}$ for all $i$. A real polynomial $\sum_{i=0}^da_it^i$ is unimodal (symmetric) if the sequence of coefficients is
unimodal (symmetric). We call a sequence $a_0,a_1,\ldots,a_d$ (or a polynomial $\sum_{i=0}^da_it^i$) strictly unimodal if,
in case $d=2m$, $$a_0<a_1<\cdots<a_m>a_{m+1}>\cdots>a_{2m},$$
and, if $d=2m+1$, 
$$a_0<a_1<\cdots<a_m=a_{m+1}>a_{m+2}>\cdots>a_{2m+1}.$$

\bigskip
There is a proof of the following result in \cite[Proposition 1]{St}.

\begin{proposition}\label{symmuni} The product of two symmetric and unimodal polynomials is symmetric and unimodal.
\end{proposition}

\noindent
{\bf Definition} For $d\ge1$, let $f_d=1+t+t^2+\cdots+t^{d-1}+t^d$. 

\begin{proposition}\label{ff}
Let $m\ge2$ be any integer, and let $1\le d_1\le d_2\cdots\le d_m$. Set $d=\sum_{i=1}^md_i$. If $f=\prod_{i=1}^mf_{d_i}=\sum_{i=0}^da_it^i$, 
then $f$ is either strictly unimodal or $a_i$ is strictly increasing until $i=\sum_{i=1}^{m-1}d_i$, then constant until $i=d_m$, then strictly decreasing until $i=d$.
\end{proposition}

\begin{proof}
Proposition \ref{symmuni} gives that $f$ is symmetric.
If $m=2$, then $a_i$ strictly increases until $i=d_1$, is then constant until $i=d_2$, then strictly decreases until $i=d_1+d_2$. We set $f_{d_1}f_{d_2}=\sum_{i=0}^{d_1+d_2}a_it^i$,
and $f_{d_1}f_{d_2}f_{d_3}=\sum_{i=0}^{d_1+d_2+d_3}b_it^i$. Then $b_i=\sum_{j=i-d_3}^ia_j$. (We let $a_j=0$ if $j<0$ or $j>d_1+d_2$.)
For $m=3$ we have two options. Either $d_2-d_1<d_3<d_1+d_2$ or $d_3\ge d_1+d_2.$ In the first case, using that $a_i$ strictly increases, then is constant, then strictly decreases,
and is symmetric, we get that $b_i=\sum_{j=i-d_3}^ia_j$ has maximum when
the interval $[i-d_3,i]$ lies as symmetrically as it can over the interval $[d_1,d_2]$, so we have one or two maxima depending on parities of $d$ and $d_3$, so $f_1f_2f_3$ is strictly unimodal.
In the second case, $d_3\ge d_1+d_2$, we get that $b_i$ increases strictly until $i=d_1+d_2$, then is constant until $i=d_3$ then decreases strictly until $i=d_1+d_2+d_3$.
We then make induction over $m$. Suppose that $f=\prod_{i=1}^mf_{d_i}$ is either strictly unimodal or that the coefficients strictly increases until $i=\sum_{j=1}^{m-1}d_j$, then is constant
until $i=d_m$, then strictly decreases until $i=\sum_{j=1}^md_j$. First suppose that $d_{m+1}<\sum_{i=1}^md_i$. If $f$ is strictly unimodal it strictly increases until one maximum 
or two consecutive maxima, then strictly decreases. The interval $I=[i-d_{m+1},i]$ has at least two elements, so $b_i=\sum_{a_i\in I} a_i$ is maximal when $I$ lies as symmetric as
possible over the maximal part of $f$. Thus also here $\prod_{i=1}^{m+1}f_{d_i}$ is strictly unimodal. Let still $d_{m+1}<\sum_{i=1}^md_i$, but now suppose
 that the coefficients $a_i$ strictly increases until $i=\sum_{j=1}^{m-1}d_j$, then is constant until $i=d_m$, then strictly decreases until $i=\sum_{j=1}^md_j$. Since  
 $d_{m+1}>d_m-\sum_{i=1}^{m-1}d_i$, we still get that $b_i=\sum_{a_j\in I} a_j$ is maximal when $I$ lies as symmetric as
possible over the maximal part of $f$, so $\prod_{i=1}^{m+1}f_{d_i}$ is strictly unimodal. Finally suppose that $d_{m+1}\ge\sum_{i=1}^md_i$.
Then $b_i=\sum_{j=i-d_{m+1}}^ia_j$ strictly increases until $i=d$, then is constant until $i=d_{m+1}$, then strictly decreases until $i=\sum_{j=1}^{m+1}d_j$.
\end{proof}

\begin{theorem}\cite[Reid,~Roberts, Roitman]{R-R-R}\label{RRR}
Let $P(t)=\sum_{i=0}^{\deg(P(t))}a_it^i$ be the Hilbert series of the graded complete intersection $k[x_1,\ldots,x_n]/(f_1,\ldots,f_n)$,
$\deg(f_i)=d_i+1$, and suppose that $d_1\le d_2\le\cdots\le d_n$. 
Then $P$ is symmetric, and either strictly unimodal or the coefficient $a_i$ strictly increases until $i=\sum_{i=1}^{n-1}d_i$, then are constant until $i=d_n$,
then strictly decreases.
\end{theorem}

\begin{proof}
The Hilbert series is $\prod_{i=1}^n(1-t^{d_i+1})/(1-t)^n$. We have $(1-t^{d+1})/(1-t)=
1+t+t^2+\cdots+t^{d-1}+t^d=f_d$. Thus the Hilbert series is $\prod_{i=1}^nf_{d_i}$. The result follows from Proposition \ref{ff}.
\end{proof}

\begin{corollary}
If $d_n+1\le\sum_{i=1}^{n-1}d_i$, then $\prod_{i=1}^nf_{d_i}$ is strictly unimodal.
In particular $f_d^n$ is strictly unimodal for any $d$ and $n\ge2$. If $d_n+1>\sum_{i=1}^{n-1}d_i$, there are $d_n-\sum_{i=1}^{n-1}d_i+1$ maximal values
for the coefficients of $\prod_{i=1}^nf_{d_i}$.
\end{corollary}

\begin{proof}
We know that $\prod_{i=1}^nf_{d_i}$ is strictly unimodal if $d_n\le\sum_{i=1}^{n-1}d_i+1$. If $d_n>\sum_{i=1}^{n-1}d_i$ we get maximal value
for $i=e,e+1,e+2,\ldots,d_n$, where $e=\sum_{i=1}^{n-1}d_i$.
\end{proof}

\section{The Betti numbers of $n+1$ generic forms in $n$ variables}
We will always assume that $g_1,\ldots,g_{n+1}$ minimally generate $(g_1,\ldots,g_{n+1})$.

\begin{lemma}
Suppose $g_i$, $i=1,\ldots, n+1$, are generic forms in $\mathbb C[x_1,\ldots,x_n]$, $\deg(g_i)=e_i$, $1\le i\le n+1$.
If $g_1,\ldots,g_{n+1}$
minimally generates $(g_1,\ldots,g_{n+1})$, then $e_{n+1}\le\sum_{i=1}^n(e_i-1)$.
\end{lemma}

\begin{proof}
The complete intersection $\mathbb C[x_1,\ldots,x_n]/(g_1,\ldots,g_n)$ has socle degree $\sum_{i=1}^n(e_i-1)$, so $\deg(g_{n+1})\le\sum_{i=1}^n(e_i-1)$.
\end{proof}

\medskip
If $I$ is a graded ideal in $S=k[x_1,\ldots,x_n]$, then $S/I$ has a graded minimal resolution

$$0\leftarrow S/I\leftarrow S\leftarrow \bigoplus_{j\ge2}S[-j]^{\beta_{1,j}}\leftarrow\cdots\leftarrow\bigoplus_{j\ge p+1}S[-j]^{\beta_{p,j}}\leftarrow 0$$
where $S[-j]$ means that the degrees of $S$ are shifted so that $S[-j]_d=S_{d-j}$. The numbers $\beta_{i,j}$ are called the graded Betti numbers of $S/I$. The
regularity of $S/I$ is $\max\{j-i;\beta_{i,j}\ne0\}$.
Each degree in the resolution gives a finite exact sequence of vector spaces. Using that the alternating sum of the dimensions in an exact sequence of vector spaces is 0, we get
the well known formula for the Hilbert series $S/I(t)$ of $S/I$ 
$$S/I(t)=\sum_{i,j}(-1)^i\beta_{i,j}t^j/(1-t)^N.$$

\begin{lemma}\label{det}
Let $R={\mathbb C}[x_1,\ldots,x_n]/(g_1,\ldots,g_{n+1})$, $g_i$ generic, $\deg g_i=e_i$. In case the first nonpositive coefficient of 
$\prod_{i=1}^{n+1}(1-t^{e_i})/(1-t)^n$ is equal to 0, the graded Betti numbers of $R$ are determined by the Hilbert series of $R$.
\end{lemma}

\begin{proof}
If the first non-positive coefficient is 0, then all $\beta_{i,j}$ with $j-i$ smaller than the regularity $D$, come from the Koszul complex on $g_1,\ldots,g_{n+1}$ 
according to \cite{Pa-Ri,Pa-Ri1} and \cite{Di}. If $j-i>D$
then $\beta_{i,j}=0$. The Hilbert series of $R$ equals $\sum_{i,j}(-1)^i\beta_{i,j}t^j/(1-t)^n$. Let the Hilbert series be $\sum_{i=0}^{n+D}c_it^i/(1-t)^n$.
Now $c_{i+D}=\sum_{j\ge i}(-1)^j\beta_{j,i+D}$. It follows from the first paragraph of Section 2 that all terms except $\beta_{i,i+D}$ come from the Koszul complex on $g_1,\ldots,g_{n+1}$, 
which are easy to determine. Thus we can determine $\beta_{i,i+D}$ from the Hilbert series.
\end{proof}

\medskip
In order not to have too heavy notation we let $e_i-1=E_i$. Recall that $f_d=1+t+t^2+\cdots+t^{d-1}+t^d$.

\begin{lemma}\label{hilb}
Let $g_1,\ldots,g_{n+1}$ be generic forms, $\deg(g_i)=e_i$, $i=1,\ldots,n+1$. Then the Hilbert series of $\mathbb C[x_1,\ldots,x_n]/(g_1,\ldots,g_{n+1})$
equals $[\prod_{i=1}^nf_{E_i}(1-t^{e_{n+1}})]$.
\end{lemma}

\begin{proof}
The Hilbert series is $[\prod_{i=1}^{n+1}(1-t^{e_i})/(1-t)^n]$ according to \cite{St}. We use that $(1-t^{e_i})/(1-t)=f_{E_i}$ for $i=1,\ldots,n$.
\end{proof}

\begin{lemma}\label{main}
Let $g_1,\ldots,g_{n+1}$ be generic forms, $\deg(g_i)=e_i$, $i=1,\ldots,n+1$.
Let $f=\prod_{i=1}^nf_{E_i}$.
Suppose $deg(f)=2l+1$ and $f=\sum_{i=0}^{2l+1}a_it^i$ with $a_0<a_1\cdots<a_{l-k}=\cdots=a_{l+k+1}>a_{l+k+2}>\cdots>a_{2l+1}$.
Then the first non-positive coefficient of $(1-t^{e_{n+1}})f$ is 0 if and only if $e_{n+1}\in\{2,4,\ldots,2k\}\cup\{3,5,\ldots,2l+1\}$. Assume
 $\deg(f)=2l$ and $f=\sum_{i=0}^{2l}a_it^i$ with $a_0<a_1\cdots<a_{l-k}=\cdots=a_{l+k}>a_{l+k+1}>\cdots>a_{2l}$.
Then the first non-positive coefficient of $(1-t^{e_{n+1}})f$ is 0 if and only if $e_{n+1}\in\{3,5,\ldots,2k-1\}\cup\{2,4,\ldots,2l\}$.
\end{lemma}

\begin{proof}
Suppose $\deg(f)=2l+1$. Now $(1-t^{e_{n+1}})f=\sum_{i=0}^{2l+1+e_{n+1}}(a_i-a_{i-e_{n+1}})t^i$, where $a_j=0$ if $j<0$ or if $j>2l+1$.
We use that the sequence $a_0,\ldots,a_{2l+1}$ is symmetric and unimodal.
First assume that $e_{n+1}\le 2k+1$. Then $a_i-a_{i-e_{n+1}}>0$ until $i=l-k+e_{n+1}$, when it is 0. Now suppose $e_{n+1}>2k+1$ and odd.
Then $a_i-a_{i-e_{n+1}}>0$ until $i$ and $i-e_{n+1}$ are placed symmetrically around the midpoint $l+1/2$. If $e_{n+1}>2k+1$ and even, then
$a_i-a_{i-e_{n+1}}\ne0$ for all $i$. If $\deg(f)$ is even, the proof is similar.
\end{proof}

\begin{theorem}\label{par}
Let $g_1,\ldots,g_{n+1}$ be generic forms, $\deg(g_i)=e_i$, $i=1,\ldots,n+1$. Then we can determine all Betti numbers of $\mathbb C[x_1,\ldots,x_n]/(g_1,\ldots,g_{n+1})$
if $\sum_{i=1}^{n+1}e_i-n$ is even. If $\prod_{i=1}^nf_{E_i}$ is not strictly unimodal, we can determine all Betti numbers also for a few cases when $\sum_{i=1}^{n+1}e_i-n$ is odd.
\end{theorem}

\begin{proof}
We use Lemma \ref{det}, Lemma \ref{hilb}, and Lemma \ref{main}. We can determine all Betti numbers if the first non-positive coefficient in $\prod_{i=1}^nf_{E_i}(1-t^{e_n+1})$
is 0. This happens when $\sum_{i=1}E_i=\sum_{i=1}^ne_i-n$ and $e_{n+1}$ have the same parity, so when $\sum_{i=1}^{n+1}e_i-n$ is even. According to Lemma \ref{main}
the first non-positive coefficient is 0 also in some other cases, if $\prod_{i=1}^nf_{E_i}$ is not strictly unimodal.
\end{proof}

\begin{corollary}
If $f$ is strictly unimodal, $\deg(f)=d$, then the first nonpositive coeficient of $(1-t^{e_{n+1}})f$ is 0 if and only if $d$ and $e_{n+1}$ have the same parity.
\end{corollary}

\begin{proof}
That $f$ is strictly unimodal means that $k=0$.
\end{proof}

\begin{example}
Let $R=\mathbb C[x_1,\ldots,x_8]/(g_1,\ldots,g_9)$, where the $g_i$'s are generic and $\deg g_i=2$ for $i=1,\ldots8$, $\deg g_9=4$. 
We claim that $\beta_{i,j}(R)$ is the element in row $j-i$ and column $i$ of

\bigskip
     
     $\begin{matrix}
      &0&1&2&3&4&5&6&7&8\\
      \text{total:}&1&9&84&399&918&1176&848&315&48\\
      \text{0:}&1&\text{.}&\text{.}&\text{.}&\text{.}&\text{.}&\text{.}&\text{.}&\text{.}\\
      \text{1:}&\text{.}&8&\text{.}&\text{.}&\text{.}&\text{.}&\text{.}&\text{.}&\text{.}\\
      \text{2:}&\text{.}&\text{.}&28&\text{.}&\text{.}&\text{.}&\text{.}&\text{.}&\text{.}\\
      \text{3:}&\text{.}&1&\text{.}&56&\text{.}&\text{.}&\text{.}&\text{.}&\text{.}\\
      \text{4:}&\text{.}&\text{.}&8&\text{.}&70&\text{.}&\text{.}&\text{.}&\text{.}\\
      \text{5:}&\text{.}&\text{.}&48&343&848&1176&848&315&48\\
      \end{matrix}$

\bigskip
The Hilbert series is $[(1-t^2)^8(1-t^4)/(1-t)^8]=[(1+t)^8(1-t^4)]=1+8t+28t^2+56t^3+69t^4+48t^5=(1-8t^2+27t^4-48t^6+48t^7-273t^8+848t^9-1176t^{10}+
848t^{11}-315t^{12}+48t^{13})/(1-t)^8$. The first non-positive coefficient in $[(1+t)^8(1-t^4)]$ is ${n\choose6}-{n\choose2}=0$.
Since the regularity of $R$ is 5, all $\beta_{i,j}$ with $j-i\le4$ come from the Koszul complex on $g_1,\ldots,g_9$. They
are $\beta_{1,2}=8,\beta_{1,4}=1,\beta_{2,4}=28,\beta_{2,6}=8,\beta_{3,6}=56,\beta_{4,8}=70$. From the Hilbert series we get $\beta_{2,7}=48,
\beta_{3,8}=273+70=343,\beta_{4,9}=848,\beta_{5,10}=1176,\beta_{6,11}=315,\beta_{7,12}=48$.
\end{example}

The Hilbert series $R(t)$ of a graded algebra $R=\mathbb C[x_1,\ldots,x_n]/I$ is determined from the graded Betti numbers, $R(t)=\sum_{i,j}(-1)^i\beta_{i,j}(R)t^j/(1-t)^n$.
In general it is possible that for some $j=j_0$ we have $\beta_{i,j_0}\ne0$ for several $i$'s, so the Betti numbers cannot be determined from the Hilbert series.
This happens in our situation e.g for $\mathbb C[x_1,x_2,x_3]/(g_1,\ldots,g_4)$ if $\deg(g_1)=\deg(g_2)=2$, and $\deg(g_3)=\deg(g_4)=4$. Then both $\beta_{1,4}$
and $\beta_{2,4}$ are nonzero. It is shown in
\cite{Mi-Mi} that this phenomenon can happen also when it cannot be explained by the Koszul complex on $(g_1,\ldots,g_{n+1})$. They only study when this happens in consecutive
homological degrees, and call the corresponding Betti numbers "ghost terms". We now show that we can have $\beta_{i,j_0}\ne0$ for two different $i$'s many times when this
cannot be explained by the Koszul complex.

\begin{theorem}
Let $I=(g_1,\ldots,g_{2k+1})\subseteq S=\mathbb C[x_1,\ldots,x_{2k}]$, $g_i$ generic for all $i$, $\deg(g_i)=2$, $i=1,2,\ldots,2k$, $\deg(g_{2k+1})=2k$.
Let the Hilbert series of $S/I$ be $\sum_{i=0}^{4k-1}c_it^i/(1-t)^{2k}$. We have $c_{2k}=-\beta_{1,2k}+(-1)^k\beta_{k,2k}$,
$c_{2k+2}=-\beta_{3,2k+2}+(-1)^{k+1}\beta_{k+1,2k+2},\ldots, c_{4k-4}=-\beta_{2k-3,4k-4}+\beta_{2k-2,4k-4}$. These do not come
from Koszul relations on $g_1,\ldots,g_{n+1}$, because there is no element of degree $(2i+1,2(i+k))$ , $1\le i\le k-2$, in this complex. 
\end{theorem}

\begin{proof}
The Hilbert series of $S/I$ is $(1+t)^{2k}-t^{2k}=((1-t^2)^{2k}-t^{2k}(1-t)^{2k})/(1-t)^{2k}$. The regularity is $2k-1$. Then $\beta_{i,j}$
with $j-i<2k-1$
comes from the Koszul complex on $g_1,\ldots,g_{2k+1}$. They are $\beta_{i,2i}={2k\choose i}$, $i=1,2,\ldots,k-2$. Let the Hilbert series be
$\sum_{i=0}^{4k-1}c_it^i/(1-t)^{2k}$. If $i\le 2k-1$ then $c_{2i}=(-1)^i{2k\choose i}$ and $c_{2i+1}=0$. If $2k\le i\le4k-1$, then 
$c_{2i}=(-1)^i{2k\choose i}-{2k\choose 2i-2k}$ and $c_{2i+1}={2k\choose 2i+1-2k}$. We have $c_{2k}=-\beta_{1,2k}+(-1)^k\beta_{k,2k}$,
$c_{2k+2}=-\beta_{3,2k+2}+(-1)^{k+1}\beta_{k+1,2k+2},\ldots, c_{4k-4}=-\beta_{2k-3,4k-4}+\beta_{2k-2,4k-4}$. The degrees $(3,2k+2),
(5,2k+4),\ldots,(2k-3,4k-4)$ do not occur in Koszul complex on $g_1,\ldots,g_{n+1}$.
\end{proof}

\begin{example}
If we consider  $\mathbb C[x_1,\ldots,x_8]/(g_1,\ldots,g_9)$, $g_i$ generic, $\deg(g_i)=2$ for $i=1,\ldots,8$, $\deg(g_9)=8$, we
get the Hilbert series
 $$\frac{1-8\,T^{2}+28\,T^{4}-56\,T^{6}+69\,T^{8}+8\,T^{9}-84\,T^{10}+56\,T
      ^{11}-42\,T^{12}+56\,T^{13}-36\,T^{14}+8\,T^{15}}{\left({1-T}\right)^{8}}$$
   The Betti numbers are:   
      
$\begin{matrix}
      &0&1&2&3&4&5&6&7&8\\\text{total:}&1&9&36&84&126&126&84&36&8\\\text{0:}&1
      
      &\text{.}&\text{.}&\text{.}&\text{.}&\text{.}&\text{.}&\text{.}&\text
      {.}\\\text{1:}&\text{.}&8&\text{.}&\text{.}&\text{.}&\text{.}&\text{.}&
      \text{.}&\text{.}\\\text{2:}&\text{.}&\text{.}&28&\text{.}&\text{.}&\text
      {.}&\text{.}&\text{.}&\text{.}\\\text{3:}&\text{.}&\text{.}&\text{.}&56&
      \text{.}&\text{.}&\text{.}&\text{.}&\text{.}\\\text{4:}&\text{.}&\text
      {.}&\text{.}&\text{.}&70&\text{.}&\text{.}&\text{.}&\text{.}\\\text{5:}&
      \text{.}&\text{.}&\text{.}&\text{.}&\text{.}&56&\text{.}&\text{.}&\text
      {.}\\\text{6:}&\text{.}&\text{.}&\text{.}&\text{.}&\text{.}&\text{.}&28&
      \text{.}&\text{.}\\\text{7:}&\text{.}&1&8&28&56&70&56&36&8\\\end{matrix}
      $
      
       \bigskip
      Here we have $\beta_{3,10}$ and $\beta_{5,10}$ different from 0, and furthermore
    $\beta_{5,12}$ and $\beta_{6,12}$ different from 0. These do not come from the Koszul complex 
      on $g_1,\ldots,g_{n+1}$, because there are no elements of degree $(2k+1,8+2k)=(1,8)+(2k,2k)$, $1\le k\le2$, in that complex.
\end{example}

\section{Declarations}
There are no competing interests. No funding was recieved.


\begin{thebibliography}{99}

\bibitem{An} D.Anick, {\it Thin algebras of embedding dimension three}, J. Algebra {\bf 100} (1986), 235--239.

\bibitem{Di} C.~Diem, {\it Bounded regularity}, J. Algebra {\bf 423} (2015), 1143--1160.

\bibitem{Fr} R.~Fröberg, {\it An equality for Hilbert seriesof graded algebras}, Math. Scand. {\bf 56} (1985), 117--144.


\bibitem{Mi-Mi} J.~Migliore and R.~Mirò-Roig, {\it On the minimal free resolution of $n+1$ general forms}, Trans. Amer. Math. Soc. {\bf 355} (2003), no. 1, 1--36.

\bibitem{Mi-Mi1}J.~Migliore and R.~Mirò-Roig, {\it Ideals of general forms and the ubiquity of the weak Lefschetz property}. J. Pure Appl. Algebra {\bf 182}
(2003), 79--107.

\bibitem{Pa-Ri} K.~Pardue and B.~Richert, {\it Syzygies of semi-regular sequences} Ill. J. Math. {\bf 53} (2009), 349--364.

\bibitem{Pa-Ri1}  K.~Pardue and B.~Richert, {\it Errata for Syzygies of semi-regular sequences} Ill. J. Math. {\bf 56} No. 4, 1001--1003.

\bibitem{R-R-R}L.~Reid, L.~G.~Roberts, and M.~Roitman. {\it On complete intersections and their Hilbert functions}, Canad. Math. Bull. {\bf 34}(4) (1991), 525--535.

\bibitem{St} R.~Stanley, {\it Log-concave and unimodal sequences in Algebra, Combinatorics, and Geometry}, Ann. New York Acad. of Sci.
{\bf 576} (1989), 500--535.

\bibitem{St1} R.~Stanley, {\it Weyl groups, the hard Lefschetz theorem, and the Sperner property}, SIAM J. Algebraic discrete methods {\bf 1} (1980), 168--184.
\end{thebibliography}
\end{document}